\documentclass{article}
\usepackage{fullpage}
\usepackage{amsmath}
\usepackage{amsthm}
\usepackage{amssymb}
\usepackage{url}
\usepackage{graphicx}
\usepackage{verbatim}
\usepackage{listings}
\usepackage{url}
\usepackage{graphicx}
\usepackage{tikz}
\usetikzlibrary{arrows,chains,matrix,positioning,scopes}
\usetikzlibrary{arrows,shapes,positioning}
\usetikzlibrary{decorations.markings}
\tikzstyle arrowstyle=[scale=1]
\tikzstyle directed=[postaction={decorate,decoration={markings,
    mark=at position .65 with {\arrow[arrowstyle]{stealth}}}}]
\tikzstyle reverse directed=[postaction={decorate,decoration={markings,
    mark=at position .65 with {\arrowreversed[arrowstyle]{stealth};}}}]

\theoremstyle{definition}

\begin{document}

\title{An Invitation to Polynomiography via Exponential Series}

\author{Bahman Kalantari \\
Department of Computer Science, Rutgers University, NJ\\
kalantari@cs.rutgers.edu
}
\date{}
\maketitle

\begin{abstract}
The subject of \emph{Polynomiography} deals with algorithmic visualization of polynomial equations, having many applications in STEM and art, see \cite{Kal04}-\cite{KalDCG}. Here we consider the polynomiography of the partial sums of the exponential series.  While the exponential function is taught in  standard calculus courses, it is unlikely that properties of zeros of its partial sums are considered in such courses, let alone their visualization as science or art. The Monthly  article  Zemyan \cite{Zemyan} discusses some mathematical properties of these zeros.  Here we exhibit some fractal and non-fractal {\it polynomiographs} of the partial sums while also presenting a brief introduction of the underlying concepts.  Polynomiography establishes a different kind of appreciation of the significance of polynomials in STEM, as well as in art. It helps in the teaching of various topics at diverse levels. It also leads to new discoveries on polynomials and inspires new applications.  We also present a link for the educator to get access to a demo polynomiography software together with a module that helps teach basic topics to middle and high school students, as well as undergraduates.
\end{abstract}


%


\date{}				




\thispagestyle{empty}



\section{Introduction} Ever since introducing the term \emph{polynomiography} for the visualization of polynomial equations via iteration functions, when encountering certain polynomials I have found it  tempting  to  consider the shape of their \emph{polynomiographs} in the complex plane. The word \emph{polynomiography} is a combination of the term \emph{polynomial}, first used in the 17th century, see Barbeau \cite{Bar}, and the suffix \emph{graphy}. Polynomiography grew out of my research into the subject of polynomial root-finding, an ancient and historic subject that continues to grow and finds new applications with every generation of mathematicians and scientists.
We can  create  literally hundreds of polynomiographs for a single polynomial equation, even when restricted to the same portion of the complex plane.

One familiar with the term \emph{fractal}, coined by Mandelbrot many years earlier, might think polynomiography is just another name for fractal images. This however is not valid  and underestimates the significance of the theory and algorithms that have led to polynomiography. One reason to call an image a polynomiograph rather than a fractal is because the image may exhibit no fractal behavior no matter in what part of the Euclidean plane it is generated.  It would thus not make sense to call it a fractal.  After all, fractal properties that may be inherent in some iterations may not be present everywhere in an image and an iterative method may be well behaved in certain areas of the Euclidean plane, or exhibit no  fractal pattern anywhere.

Consider for instance the polynomiograph of the polynomial $1+z+0.5z^2$, the quadratic partial sum of the exponential function under the iterations of Newton's method,  shown in Figure \ref{Fig1}, top-left.  Given a particular point $z_0$, the iterations of Newton's method generate a sequence $z_0, z_1, z_2, \dots,$ that may or may not converge to a root. The sequence is called the \emph{orbit} of $z_0$.  The \emph{basin of attraction} of a root is the set of all points in the plane whose orbit converges to that root. In the figure the upper half represents the  basin of attraction of the one root and the lower half the basin of attraction of the other root.  The points on the perpendicular bisector of the line connecting the roots, the $x$-axis,  does not belong to either one of the two basins.
There is no fractal property in this image,  no broken lines, no self-similarity.  There is chaos corresponding to orbits of the points on the bisecting line, however this set is not a fractal set.  The basins of attraction form the \emph{Fatou} set, and the bisecting line the \emph{Julia} set.   This is an example of a non-fractal polynomiograph.  On the other hand, for a polynomial of degree three with distinct roots, the corresponding Newton's polynomiograph is fractal if it contains a portion of the Julia set.  The  Julia set is the boundary of each basin of attraction of a root and is fractal but the Fatou set may contain more than the union of the basins of attractions.  There are deep results on the Julia and Fatou sets of polynomials under the iterations of Newton's method and other iteration functions, including their own iterations.  For instance,  the animation \cite{KS} gives  a 3D depiction of the dynamics of the Fatou sets for the polynomial $z^3-1$. About fractals and the amazing mathematical properties of iterations of rational functions one can consult Beardon \cite{Bea}, Devaney \cite{Dev}, Mandelbrot \cite{Man}, Milnor \cite{Milnor}, and Kalantari \cite{Kalbook} with emphasis on iteration functions for polynomial root-finding.

 A polynomiograph may exhibit fractal features, in which case one can refer to it as a \emph{fractal polynomiograph}. Even if a polynomiograph exhibits fractal behavior it is more informative to refer to it as such, rather than plain fractal. Indeed the word fractal is used in very general terms. It may refer to many different types of objects, such as 2D images coming from iterations of all kinds, 3D fractals and  even objects such as trees, clouds, mountains,  nature and the universe.  After all, just because we may refer to an image or an object as fractal it does not imply that we understand all its properties.

In \cite{Kalbook} and several other articles  I have described reasons in support of the definition of the term.  While originally polynomiography  was to represent as algorithmic visualization of a polynomial equation via a specific family of iteration functions, called the \emph{basic family} to be discussed later in the article, after many more years of experiences I would like  polynomiography to refer to visualization of polynomials in broader terms,  allowing the possibility of other iterative methods, even a mixture of iterative methods,  even 3D visualizations or visualizations that pertain to the zeros of complex polynomials, not necessarily via iterations.

Based on many educational experiences, including those of educators and students who have come to experiment with polynomiography software,  there is convincing evidence that the images  convey meaningful mathematical attributes of polynomials  and algorithmic properties that make the images interesting beyond their aesthetic beauty as art, especially to the youth and students.  There have been attempts by educators to popularize fractals in education and to teach some basic properties, for instance at high schools.  However, introducing fractals in a very general setting could be confusing to the youth. After all, the underlying mathematics of fractals and iterations is sophisticated.  On the other hand,  since  solving quadratic equations is common knowledge in middle and high school, students can connect to polynomiography in an easy fashion, turning polynomials of any degree into fun objects to deal with.  Studying the underlying theory of polynomiography also makes it possible to teach and learn about fractals. Polynomiography helps motivate the teaching of fractals and related material at the K-12 level and beyond in a constructive manner, also  connecting geometry and algebra. The present article is an attempt to demonstrate the beauty of a well known class of polynomials, seldom considered as complex polynomials in the manner presented here.

What makes visualization of a polynomial equation interesting, even when all its coefficients are real numbers, is to view its domain not as the real line, but the Euclidean plane. The simplest case of polynomiography is the visualization of the basins of attraction of Newton's method when applied to a quadratic equation, e.g. $z^2-1$, historically  considered by Cayley \cite{cay} in 1879. Except for a shift, its polynomiograph is identical with the one shown in Figure \ref{Fig1}, top-left image.   It is not difficult to mathematically prove this property of the basins of attraction without computer visualization. On the other hand, the analysis of basins of attraction of $z^3-1$  under Newton's method  is very complex.  It is well known that the resulting image is fractal, a  fractal polynomiograph.  With the advent of computers it became more plausible to understand the shape of basins of attraction and apparently the first person who tried this and saw the surprising fractal behavior was the mathematician John Hubbard, see Gleick \cite{Gl}.  While polynomiographs of Newton's method for quadratic are not fractal, one can easily modify Newton's method to get fractal polynomiographs.  This modified method is a \emph{parametrized Newton's method}, described in more generality in the next section.

The subject of fractals is significant, vastly rich and very beautiful. Mandelbrot not only coined the term fractal but undoubtedly played an enormous role in bringing them into view and this in turn has resulted in many theoretical advancements and visualizations, including new kinds of algorithmic mathematical art.  Polynomiography does overlap with fractals in many way. However, it is not a subset of fractals, not in theory, nor in practice, nor in terms of its images as art or otherwise.
I believe that polynomiography can play an important role in the teaching of fractals and dynamical systems at various levels. To support this point, in numerous personal experiences that include presentations  to hundreds of middle and high school students,  lectures during a first-year seminars or formal courses,  only a very small percentage of students have ever heard of the term ``fractal.''  Even those who had familiarity with fractals could only identify them as visual images that represent self-similarity.  Even at universities, topics on fractals and dynamical systems are only offered as graduate level courses.  For those who wish to teach or learn basic concepts from the theory of fractals and dynamical systems polynomiography can provide a  powerful bridge  into these subjects areas, as well as many others. Polnomiography appeals to students because they can connect it with a task they have leaned early on, namely solving a polynomial equation. This in particular makes polynomiography effective for introducing it to K-16 students at elementary or advanced levels.

As an example consider introducing the Mandelbrot set to middle and high students beyond just showing the aesthetic  beauty of the set. We must first introduce them to Julia sets resulting from  the iterations of a quadratic function.  However, these iterations attempt to find fixed points but not roots. The notion of fixed points, while implicit in Newton's method, typically is not taught in K-12, and hence is unfamiliar to students.  In order to popularize fractals, first  the notion of fixed points and fixed point iterations must be introduced. Having introduced these, then we can consider  the task of approximating the roots of a polynomial $p(z)$ as that of finding the fixed points of the polynomial, $q(z)=p(z)+z$.   A fixed point $\theta$ of $q(z)$ is  \emph{attractive, repulsive} and \emph{indifferent} if the modulus of $p'(\theta)$ is  less than one, larger than one, or equal to one.  For a complex number $z=x+i y$, $i=\sqrt{-1}$,  its modulus is $|z|= \sqrt{x^2+y^2}$.  The fixed point iteration refers to computing $z_{k}=q(z_{k-1})$ for $k=1, 2, \dots$,  where $z_0$ is a starting \emph{seed}. If a fixed point is attractive, the iterations are guaranteed to converge to it, provided $z_0$ is close enough.   That iterations are necessary to approximate the roots even of a quadratic equation is clear to a middle schooler who knows the quadratic formula fails to provide a numeric decimal value to the solution of $x^2-2$.  Students can experience the behavior of iterations of quadratic functions via polynomiography. These experiences will demonstrate that not both fixed points of a quadratic can be attractive.  However,  Newton's method will never fail to approximate the roots because both fixed points are attractive with respect to Newton's function.  While Newton's method doesn't result in fractal polynomiographs for  any quadratic, fractal Julia sets result in polynomiography of a quadratic corresponding to a parametrized Newton's method, see Figure \ref{Fig3}. By introduction of  different values for the parameter, students quickly learn the notion of attractive and repulsive fixed points and  appreciate Newton's method for numerically solving a quadratic equation.   Students can also be introduced to the notion of open and closed sets, as well as Fatou and Julia sets.

By considering polynomiography students and teachers quickly discover the vastness of the world of polynomials and they discover new applications of them, distinct from the traditional applications considered in standard textbooks.  Among these applications one can include educational lesson plans, visual cryptography, art and much more. In  my personal experiences in the teaching of the subject of polynomiography I am often delighted to find many creative applications that students are able to think of, including those that I did not even imagine.  Why limit the applications of polynomials to the standard ones discussed in typical textbooks on algebra, calculus, or numerical analysis? Why not think of the shape of zeros of a polynomial, even if these do not come up in ordinary application? Once we think of the zeros of a polynomial, its polynomiography  becomes a relevant matter of curiosity, leading to new images, new discoveries, new applications, new questions, and new art.

Polynomiography may be a  mouthful of a word, however it is a meaningful one. Students quickly accept it.  In order to introduce polynomiography we need to consider polynomial equations over the complex plane. Since everyone is already familiar with the Cartesian coordinate system in the plane, it is easy to describe  a polynomial equation as a way to encrypt a bunch of points in the Euclidean plane. We think of the points as \emph{complex numbers}. This allows for turning points in the Euclidean plane into objects that inherit the four elementary operations on real numbers.   For a middle or high school student,  learning about elementary operations on complex numbers is a matter of minutes rather than hours. Once these operations are understood, a polynomial equation together with the \emph{fundamental theorem of algebra} is nothing more than a way to encrypt points.  Solving a polynomial equation is a game of hide-and-seek.  For a fun introduction to the fundamental theorem of algebra, see Kalantari and Torrence \cite{KT}.

In this article I present some polynomiography for the partial sums of the exponential series, familiar to every student who has come across calculus. The exponential function is considered  by some mathematicians to be the most important function in mathematics. Polynomiography for the partial sums of some analytic functions such as sine and cosine is already considered in \cite{Kalbook}.  In fact we can do polynomiography for functions that are not polynomial, and witness a visual convergence in the sense of polynomiography. The irony is that the exponential function itself has no zeros.  At the end for the educator I will provide links to demo polynomiography software and a teaching module.

\section*{The $n$-th Partial Sum of Exponential Series and The Basic Family}

The exponential function and its  $n$-th partial sum polynomial are, respectively
\begin{equation}
\exp(z)= \sum_{k=0}^\infty \frac{z^k}{k!}, \quad P_n(z)= \sum_{k=0}^n \frac{z^k}{k!}=1+z+\frac{z^2}{2!}+\cdots + \frac{z^n}{n!}.
\end{equation}
The shape of zeros of these polynomials has been studied, see Zemyan \cite{Zemyan} for a wonderful review. Many results are known, for instance, bounds on the zeros of $P_n(z)$. Some conjectures are raised on the zeros, such as the convexity of the roots are described in \cite{Zemyan}. Undoubtedly meany research questions can be stated.

The polynomiographs of the partial sums to be exhibited here are generated via the \emph{basic family} of iteration functions. For a given arbitrary polynomial $p(z)$, the basic family is an infinite  collection of iteration functions. To define the basic family even in more generally, given a complex number $\alpha$ satisfying  $\vert 1- \alpha \vert < 1$, the \emph{ parametrized basic family} is:
\begin{equation}B_{m, \alpha}(z)=z- \alpha p(z) \frac{D_{m-2}(z)} {D_{m-1}(z)}, \quad m=2,3, \dots \end{equation}
where  $D_0(z)=1$,  $D_k(z)=0$ for $k <0$, and, $D_m(z)$ satisfies the recurrence
\begin{equation}
D_m(z)= \sum_{i=1}^n (-1)^{i-1}p(z)^{i-1}\frac{p^{(i)}(z)}{i!}D_{m-i}(z).\end{equation}
The range for $\alpha$ guarantees that each root of $p(z)$ remains an attractive fixed point of $B_{m, \alpha}(z)$, see \cite{Kalbook}.  When $\alpha=1$ we denote the family by $B_m(z)$. Each member is capable of generating different polynomiography of the same polynomial.  The first two members are, the $B_2(z)$, Newton, and the $B_3(z)$, Halley, iteration functions.  This family and its variations are extensively studied in \cite{Kalbook} establishing  many fundamental properties of why they are probably the most important family of iteration functions for polynomial root-finding. In particular, for each fixed $m\geq 2$, there exists a disc centered at a root $\theta$ such that for any $z_0$ in this disc the sequence of fixed point iteration
$z_{k+1}=B_m(z_k)$, $k=0,1,\dots$, is well-defined and converges to $\theta$. When  $\theta$ is a \emph{simple
root}, i.e. $p(\theta)=0$, $p'(\theta) \not =0$,  the order of convergence is $m$. Variations of the basic family,  other than the parameterized version, are described in \cite{Kalbook}.

In contrast to using individual members of the basic family, there is a collective application, using the \emph{ basic sequence}, $\{B_m(w), m=2, \dots\}$, where $w$ is some fixed complex number, see \cite{Kalbook}. To describe the convergence of the basic sequence we need to define the notion of the \emph{Voronoi diagram} of a set of points in the Euclidean plane.  Given a set of points $\theta_1, \dots,  \theta_n$ in the Euclidean plane, the \emph{ Voronoi cell} of a particular point $\theta_i$ is the set of all points in the plane that are closer to $\theta_i$ than to any other $\theta_j$.  If $w$ lies in the Voronoi cell of a particular root $\theta$ of $p(z)$, then it can be shown that the sequence $B_m(w)$ converges to $\theta$.
For pointwise convergence  see \cite{Kalbook}, and for a proof of a stronger property, uniform convergence of the basic family, see \cite{KalDCG}. Based on this convergence property we can produce non-fractal polynomiographs that are very different from the usual fractal ones.

\section*{Polynomiographs for the Partial Sum}

\begin{figure}[h!tbp]
\begin{center}
\includegraphics[width=1.5 in]{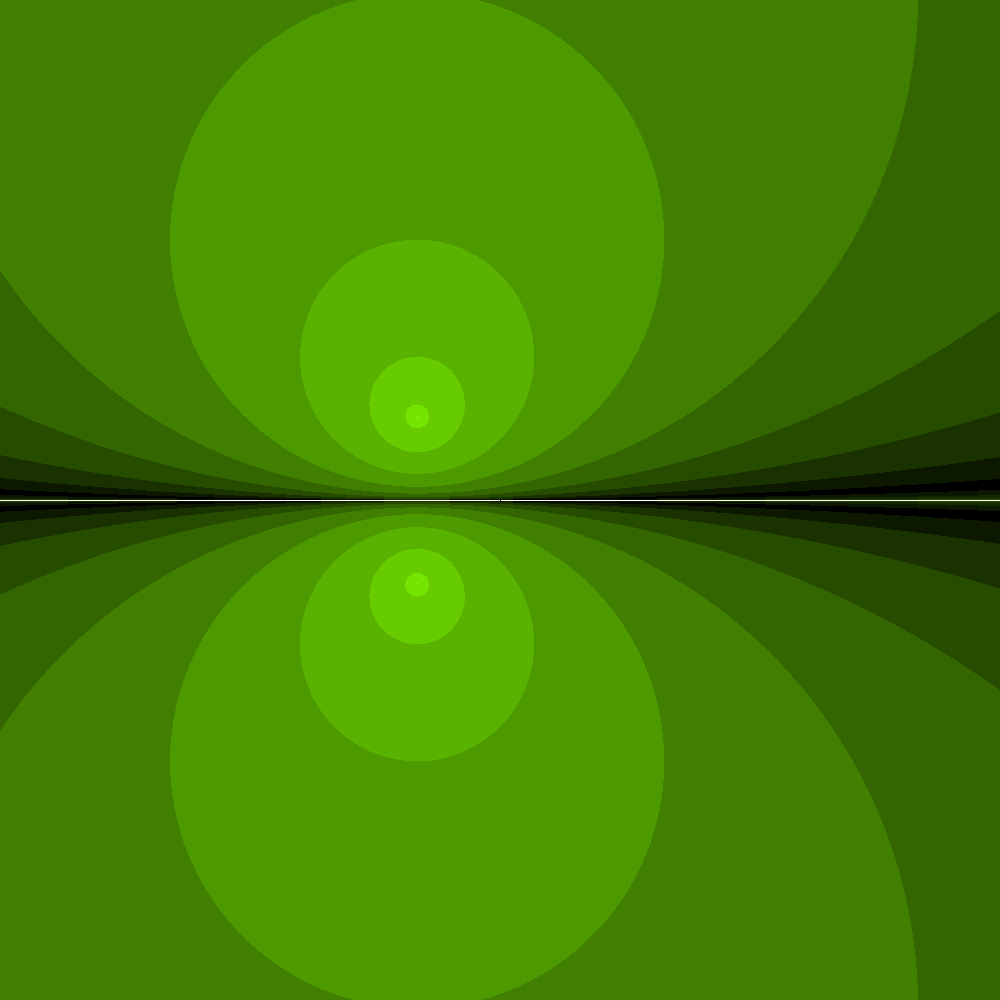} \includegraphics[width=1.5 in] {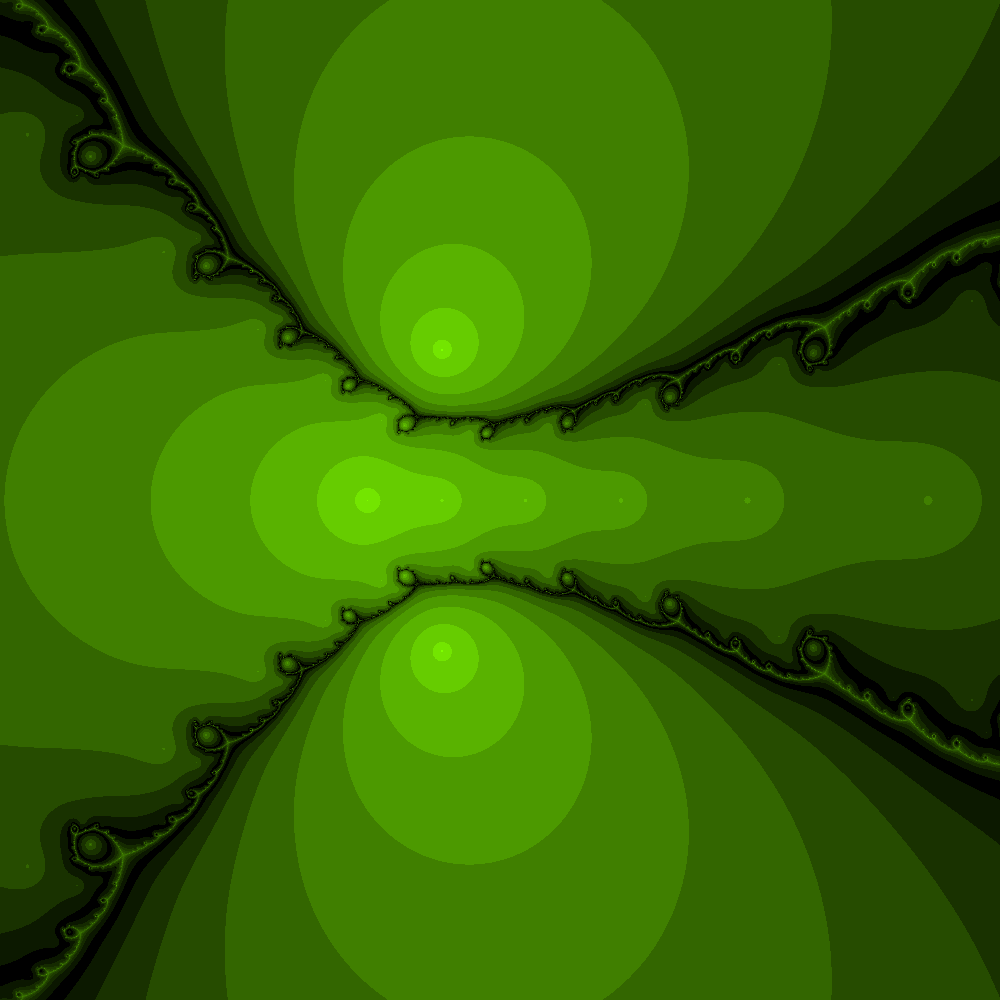} \includegraphics [width=1.5 in]{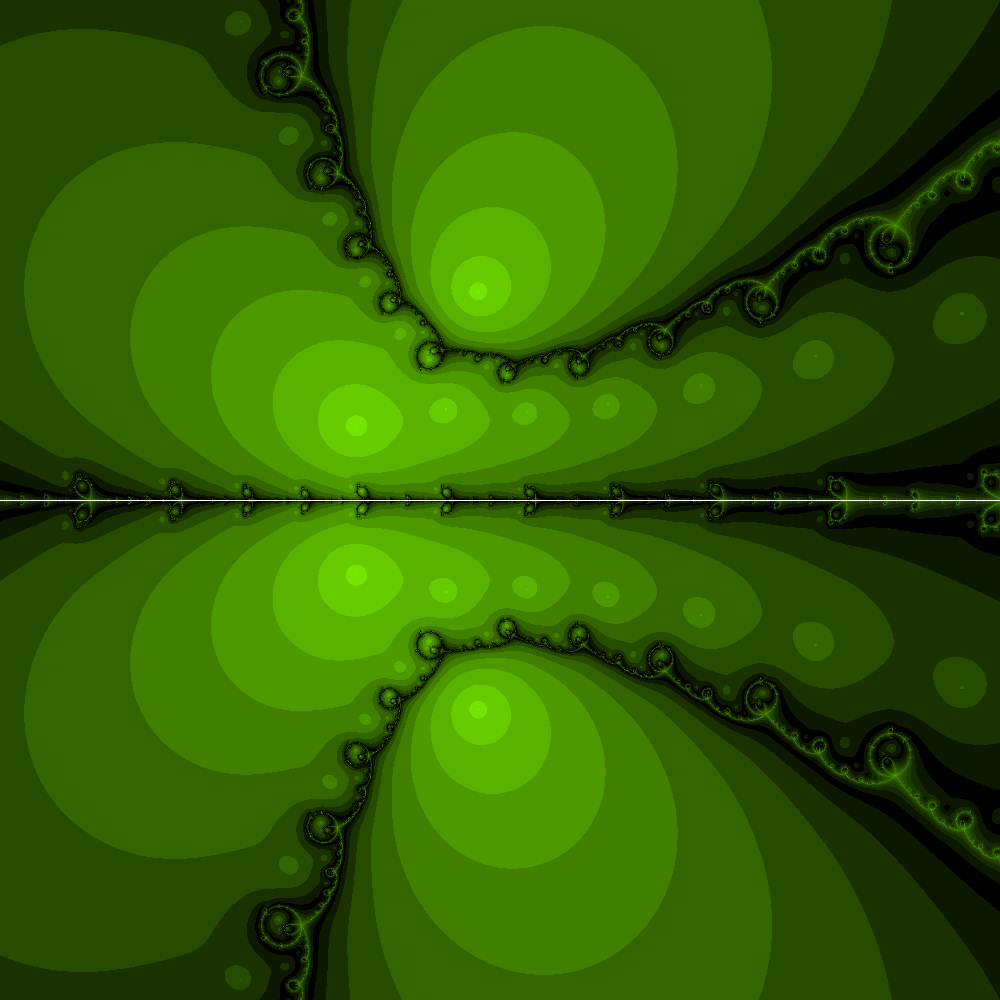}\

\vspace{0.1 cm}
\includegraphics[width=1.5 in]{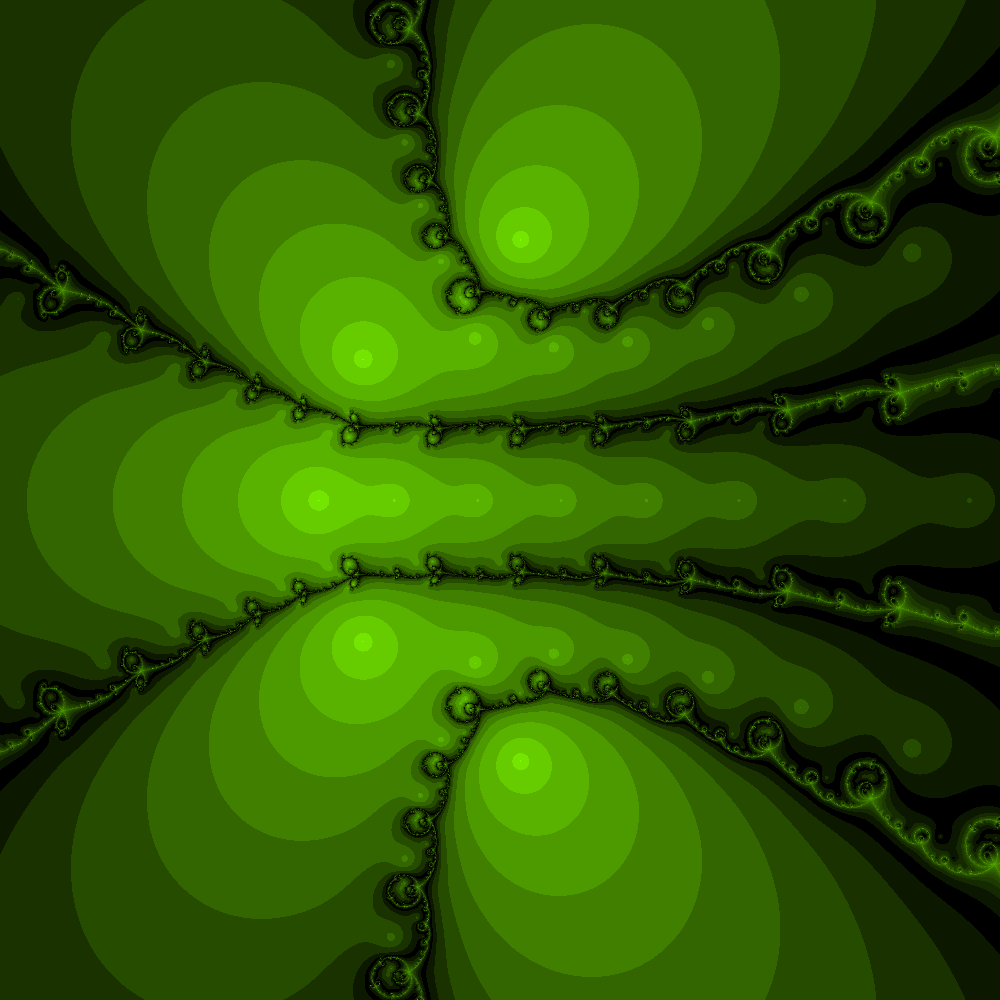} \includegraphics[width=1.5 in]{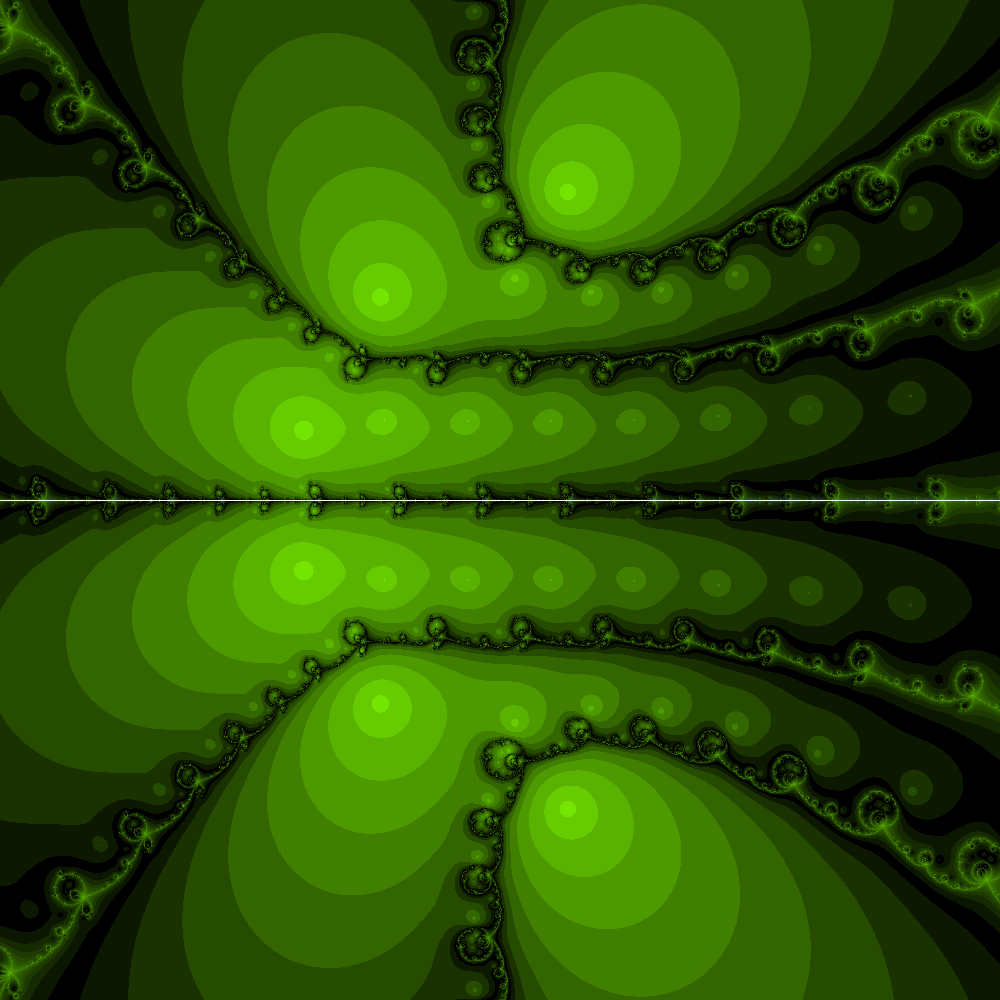} \includegraphics[width=1.5 in]{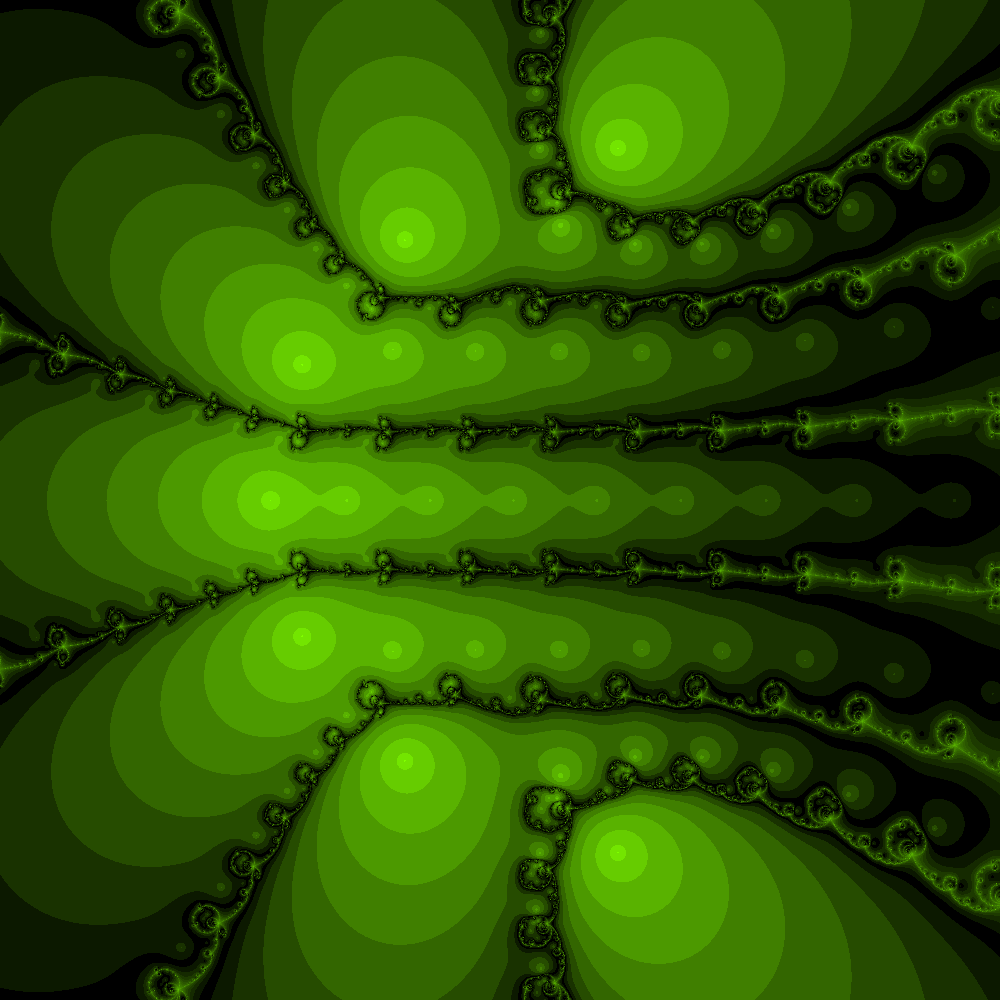}\

\vspace{0.1 cm}
\includegraphics[width=1.5 in]{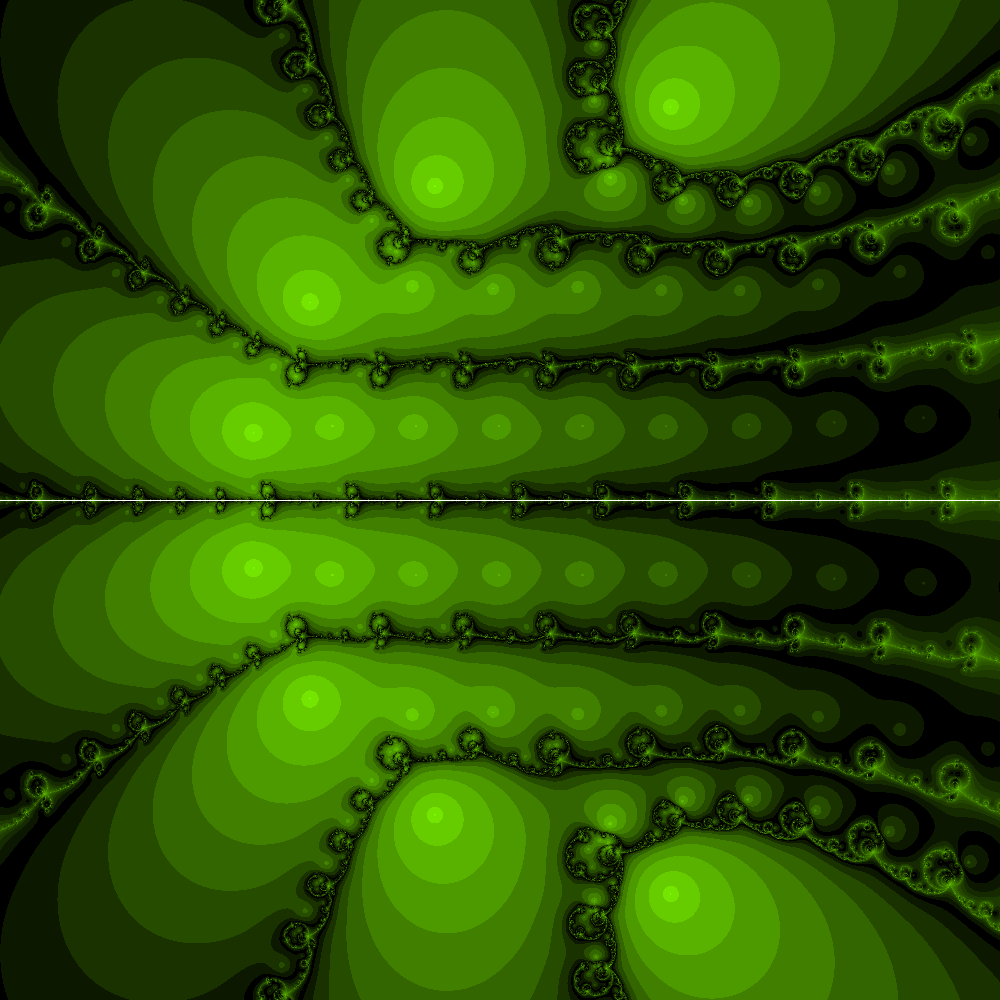} \includegraphics[width=1.5 in]{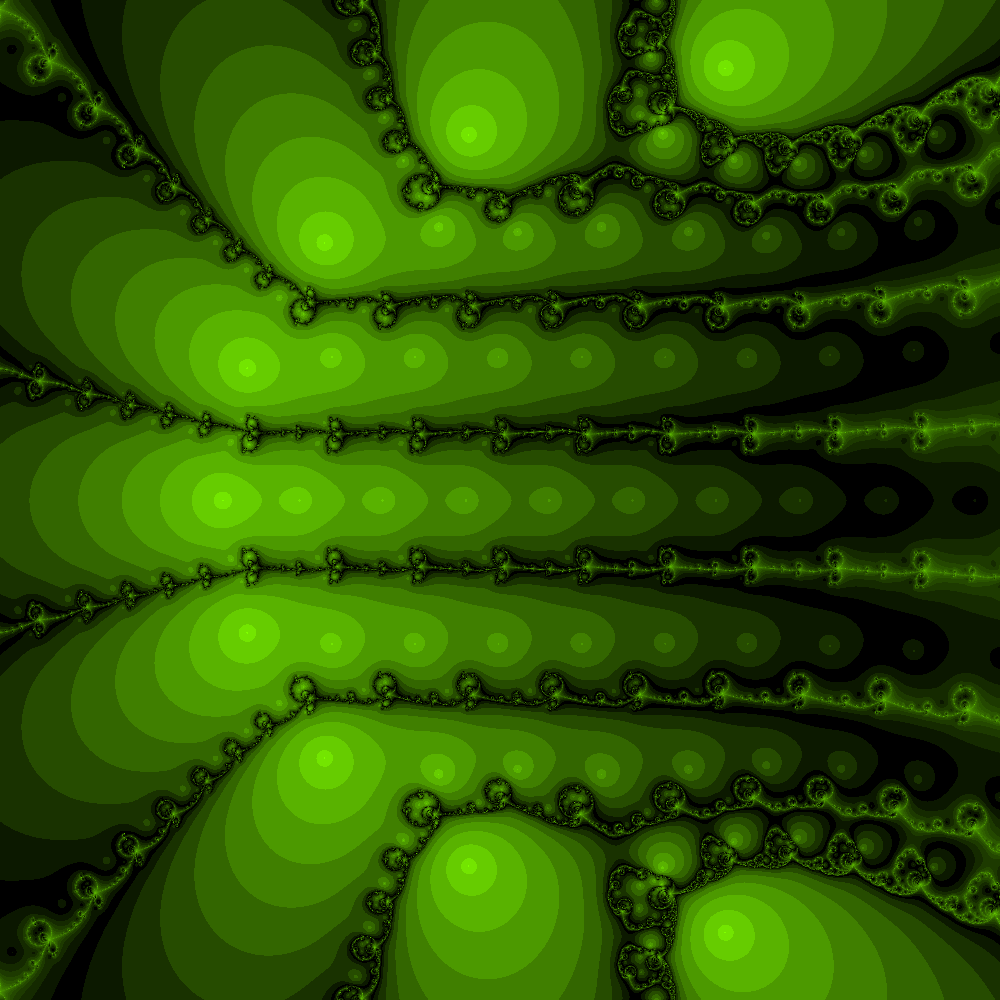} \includegraphics[width=1.5 in]{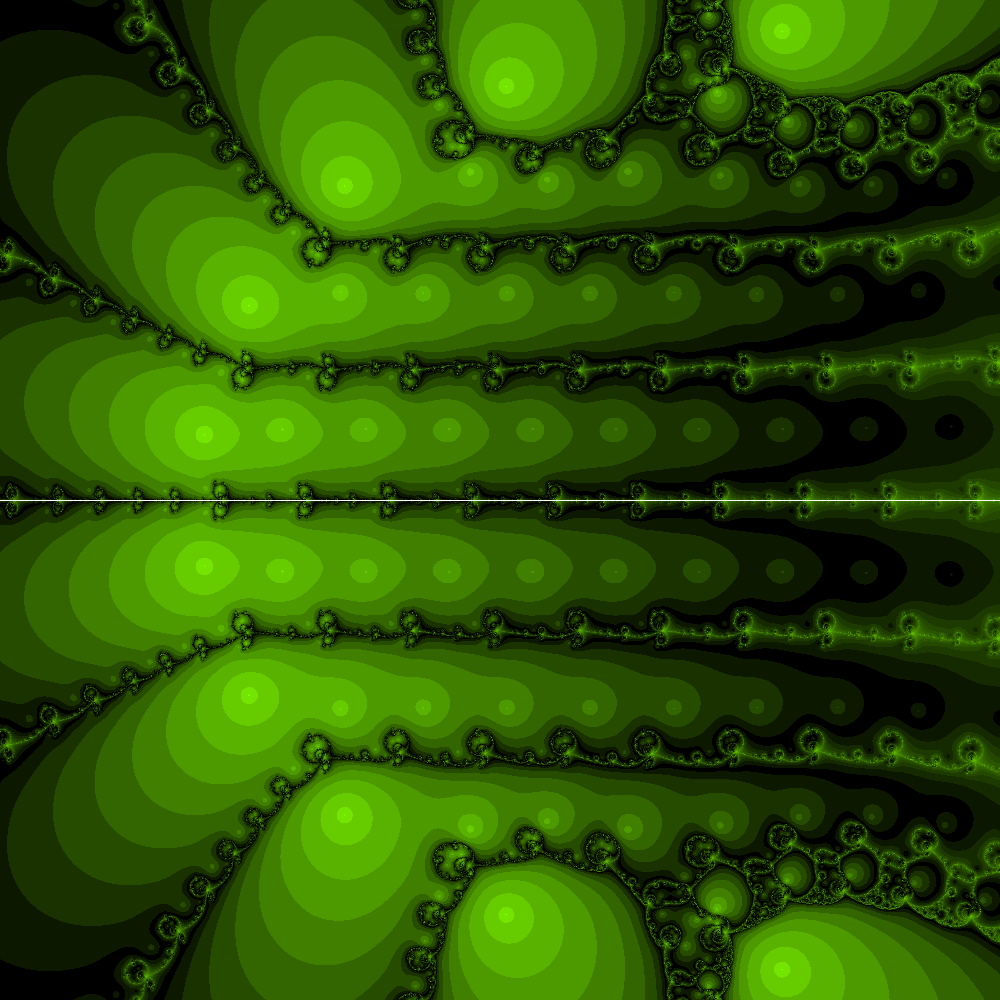}\

\vspace{0.1 cm}
\caption{Polynomiographs of $P_n(z)$, $n=2, \dots, 10$ under Newton's method.} \label{Fig1}
\end{center}
\end{figure}

\begin{figure}[h!tbp]
\begin{center}
\includegraphics[width=1.5 in]{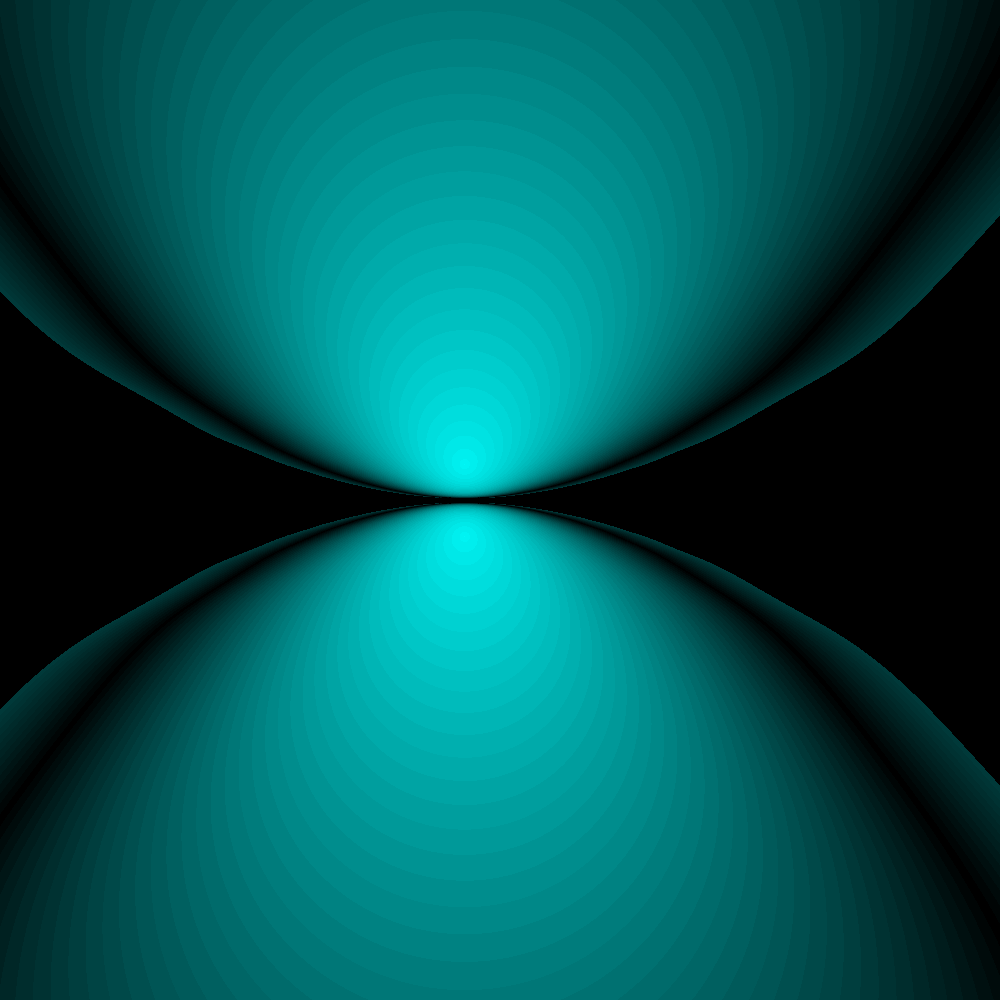} \includegraphics[width=1.5 in] {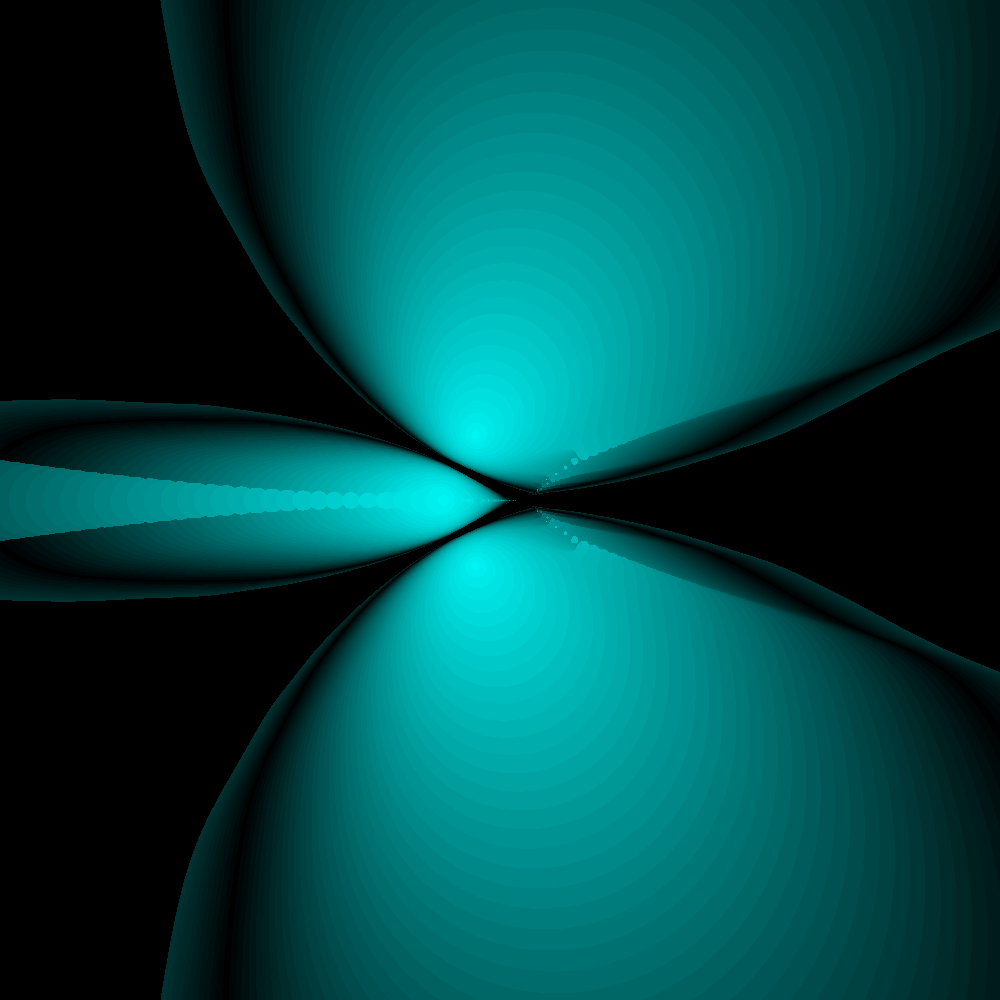} \includegraphics [width=1.5 in]{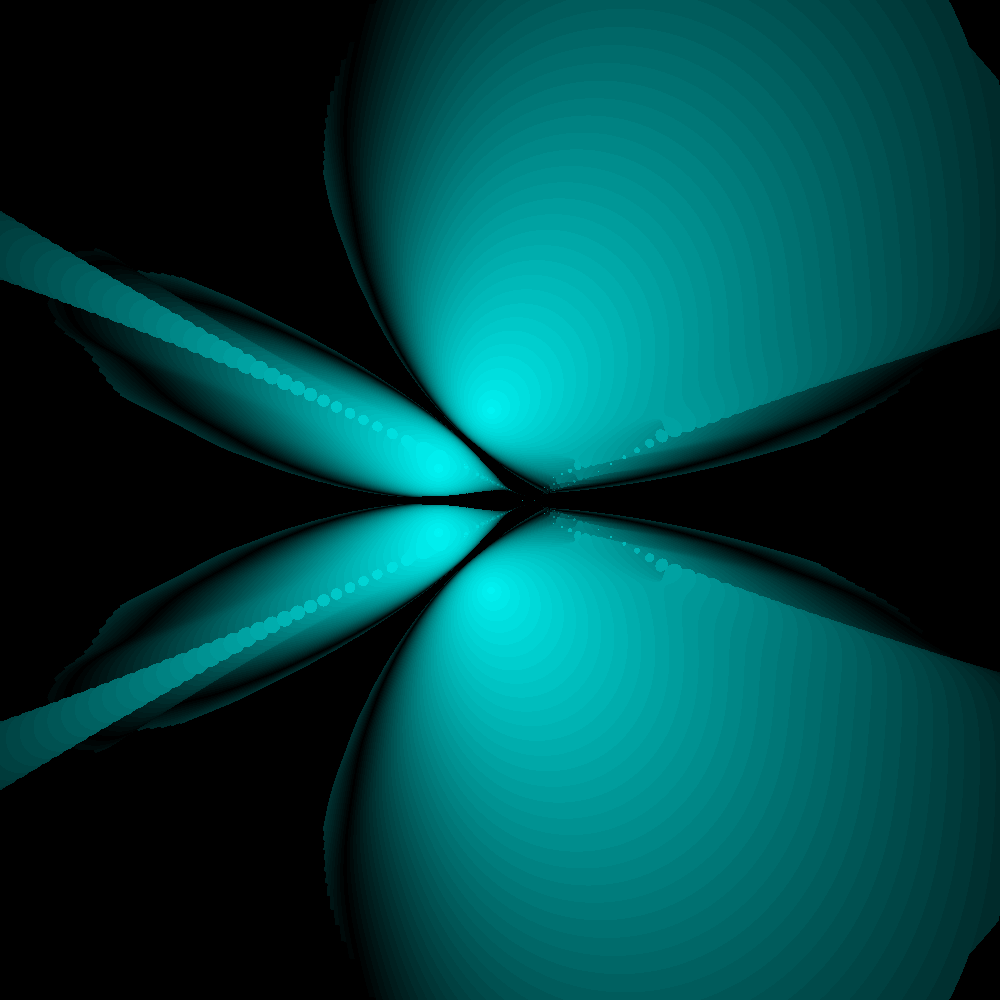}\\
\vspace{0.1 cm}
\includegraphics[width=1.5 in]{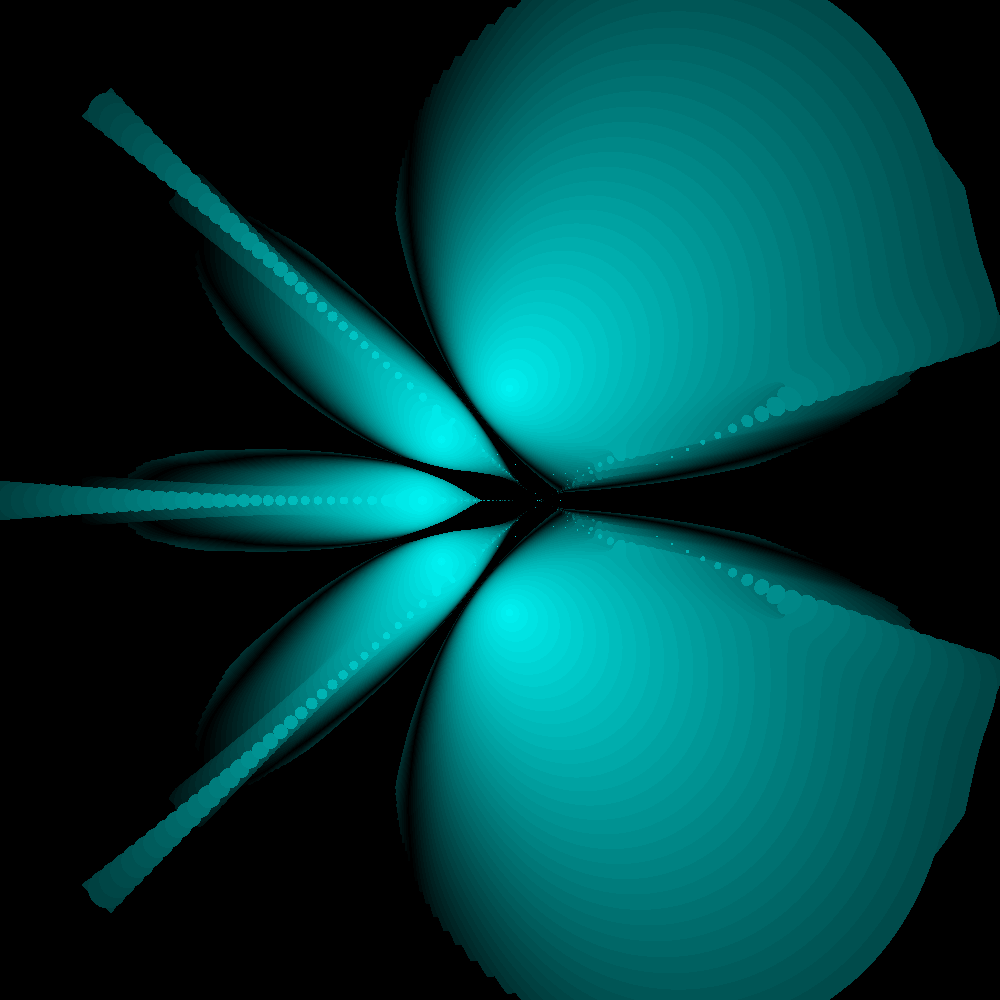} \includegraphics[width=1.5 in]{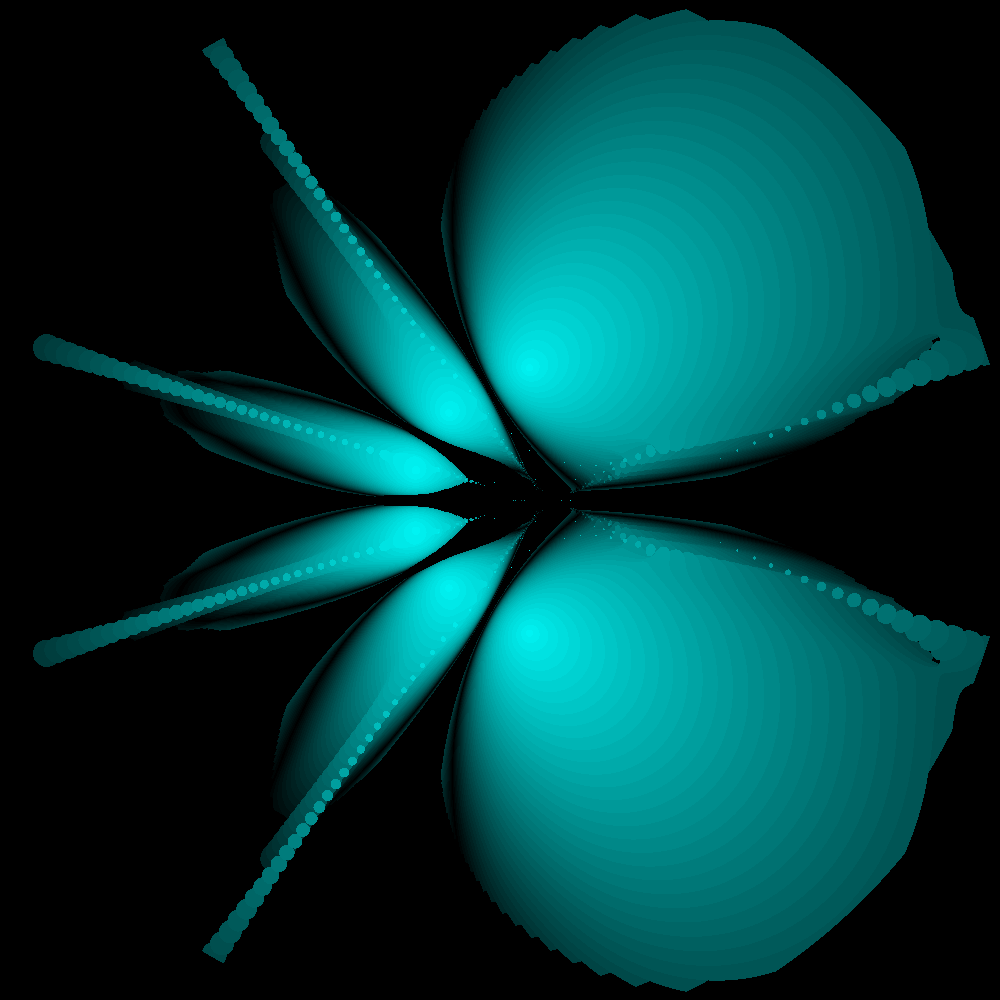} \includegraphics[width=1.5 in]{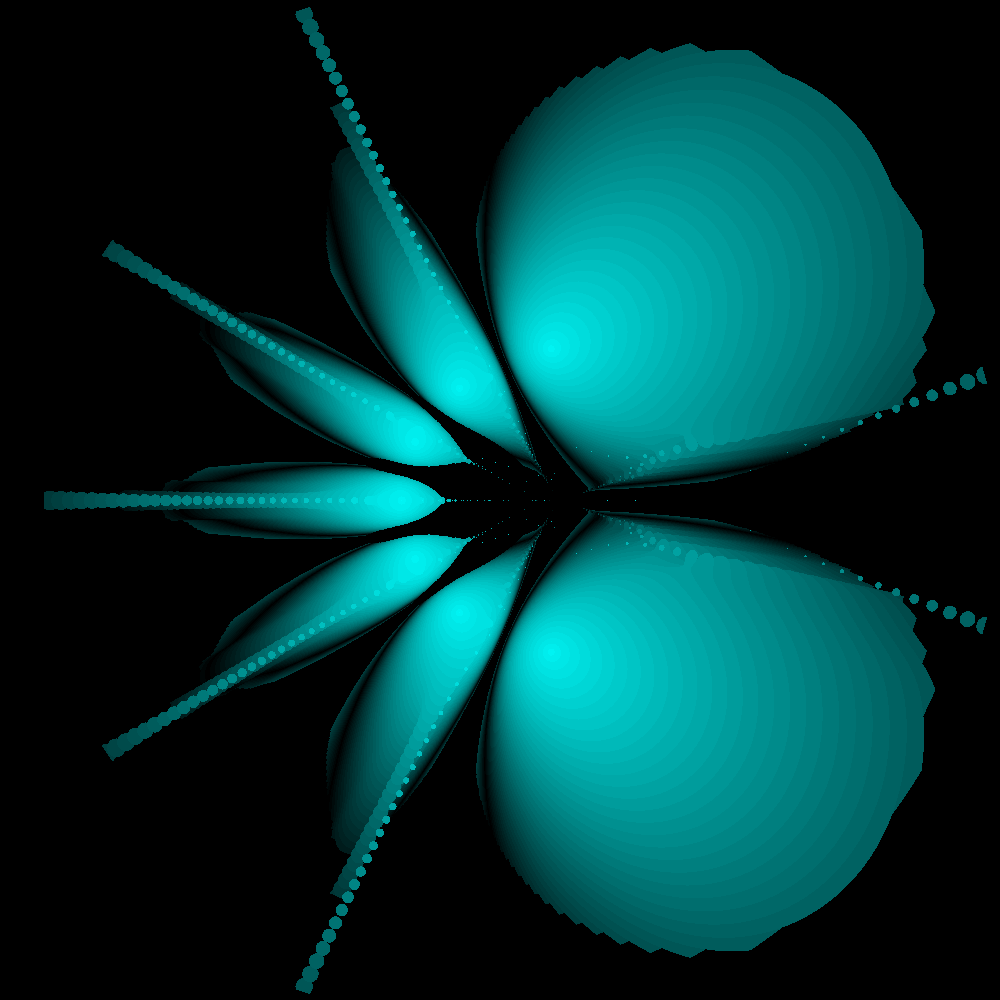}\\
\vspace{0.1 cm}
\caption{Polynomiographs of $P_n(z)$, $n=2, \dots, 7$ under point-wise convergence.} \label{Fig2}
\end{center}
\end{figure}

\begin{figure}[h!tbp]
\begin{center}
\includegraphics[width=1.5 in]{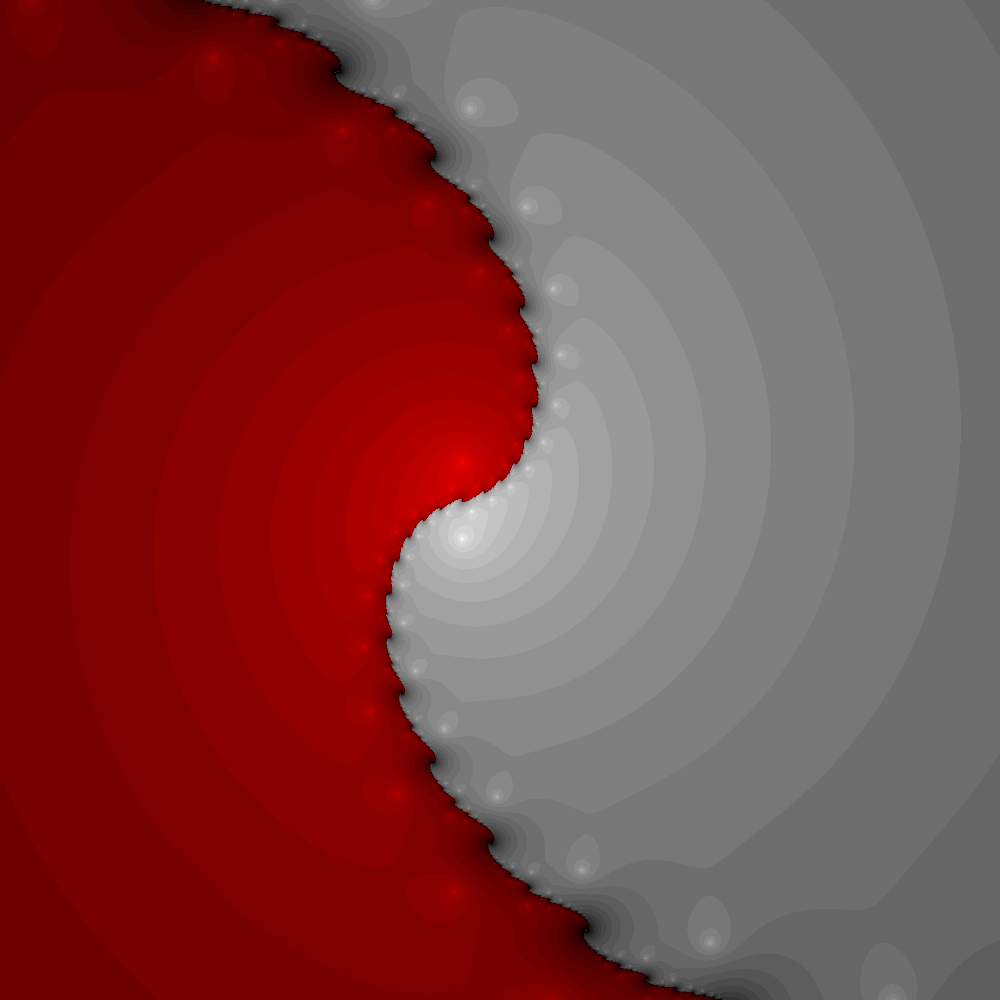} \includegraphics[width=1.5 in] {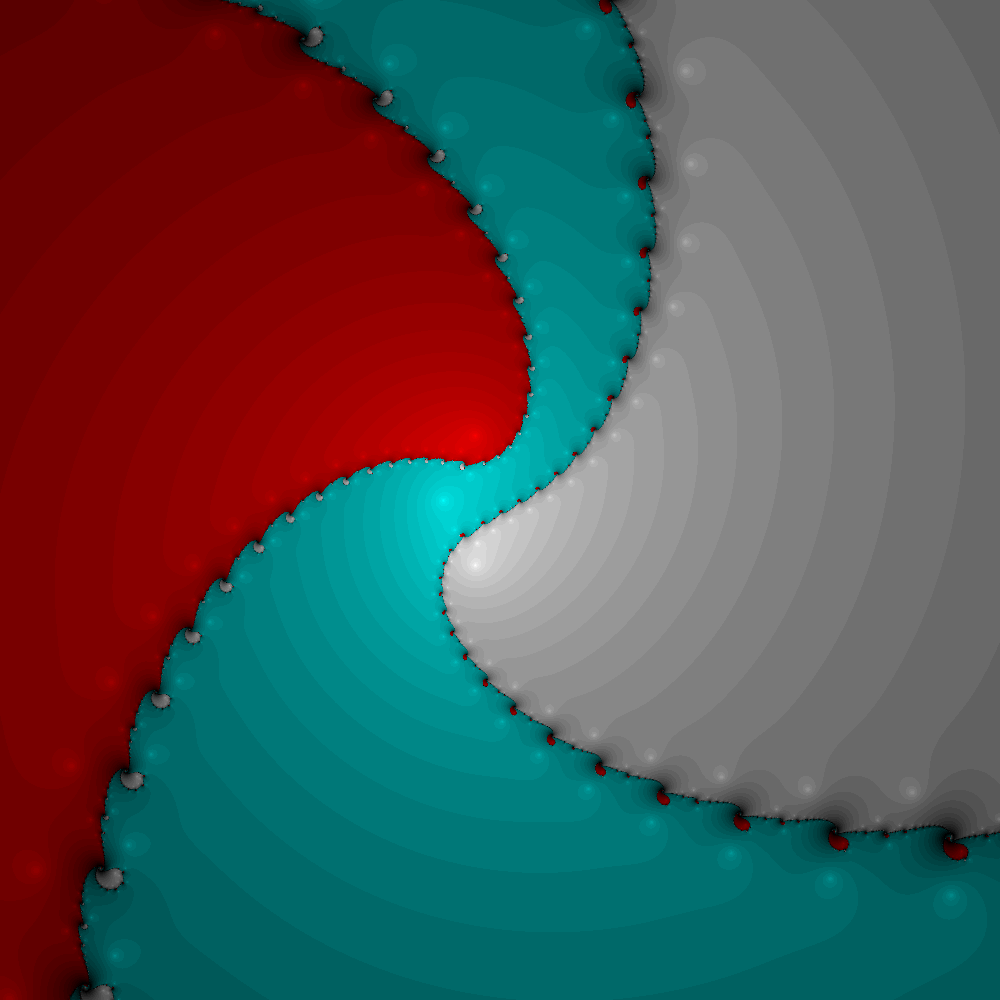} \includegraphics [width=1.5 in]{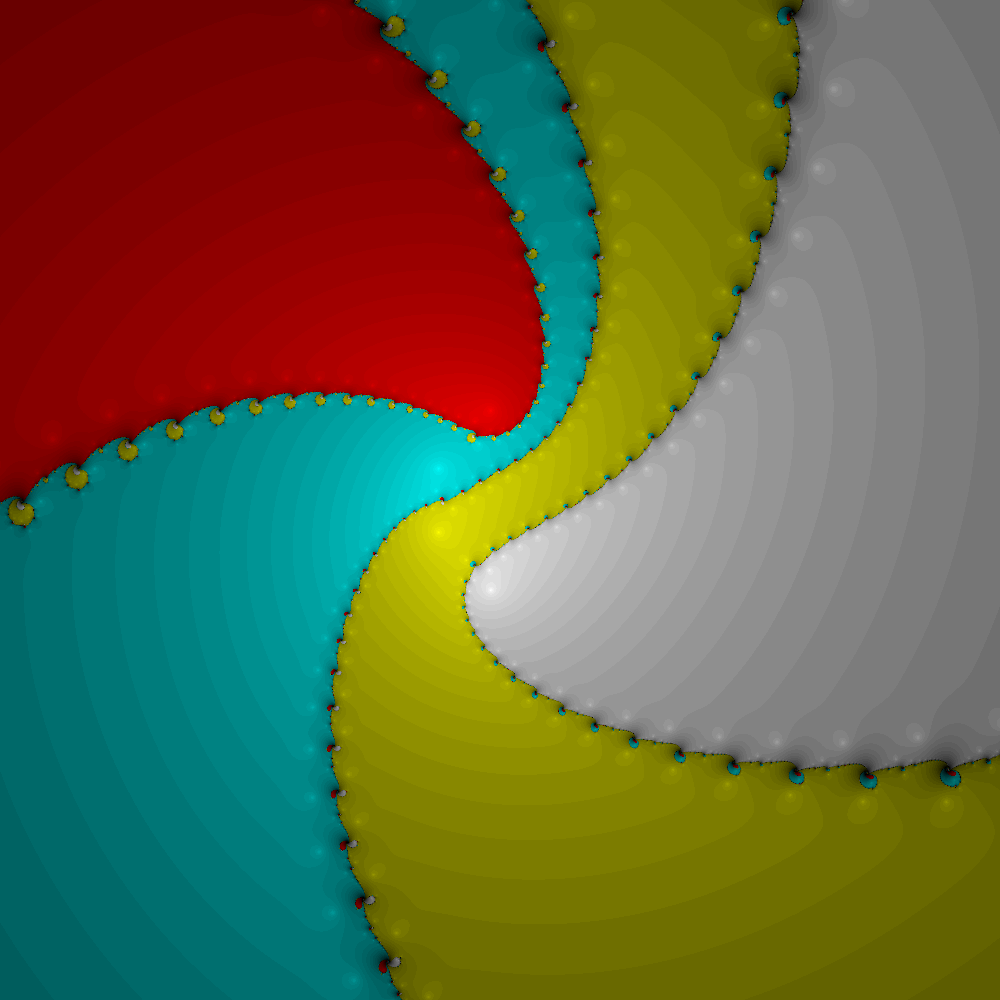}\\
\vspace{0.1 cm}
\includegraphics[width=1.5 in]{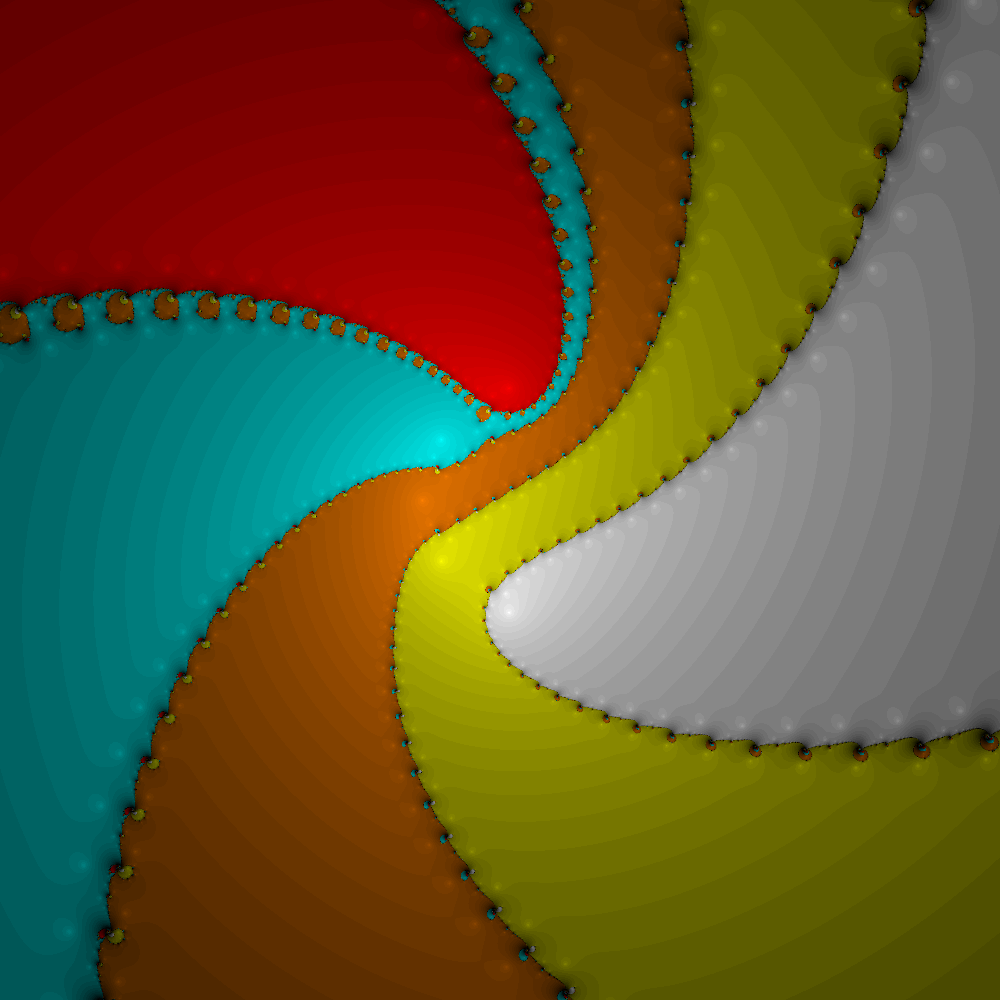} \includegraphics[width=1.5 in]{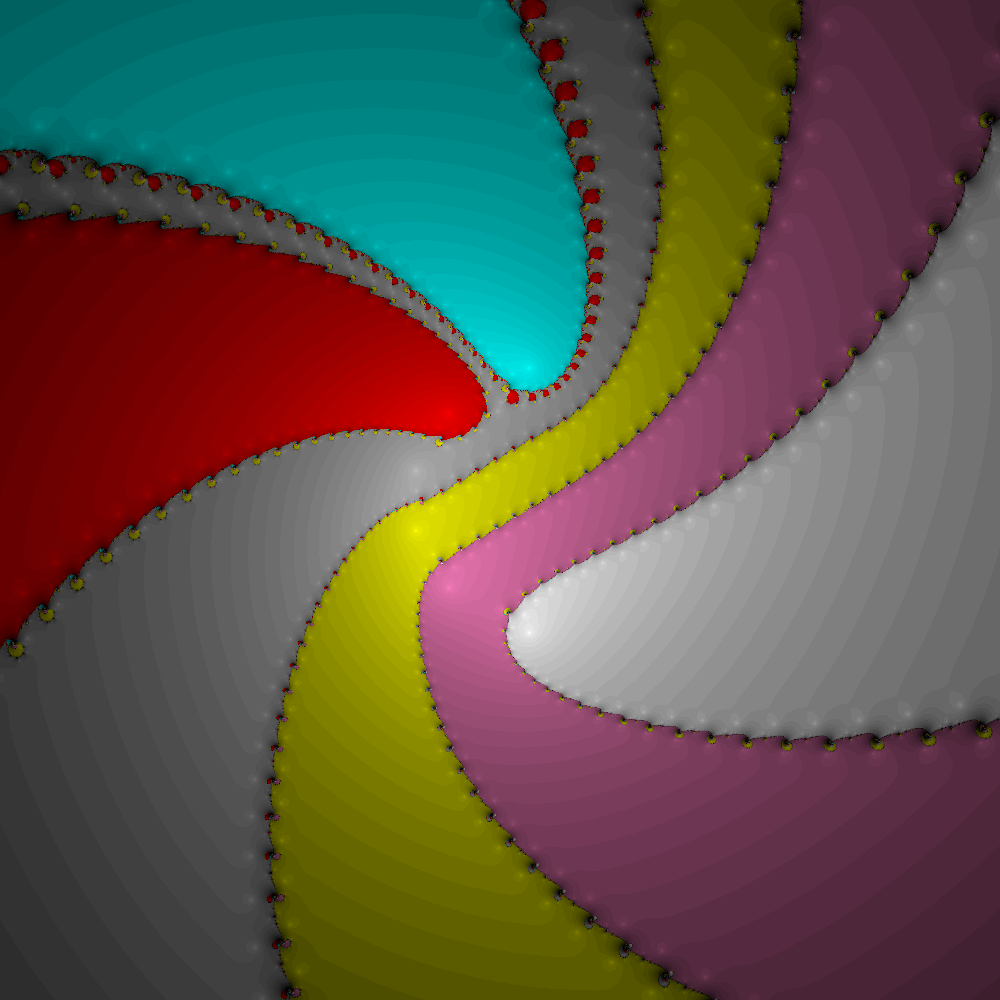} \includegraphics[width=1.5 in]{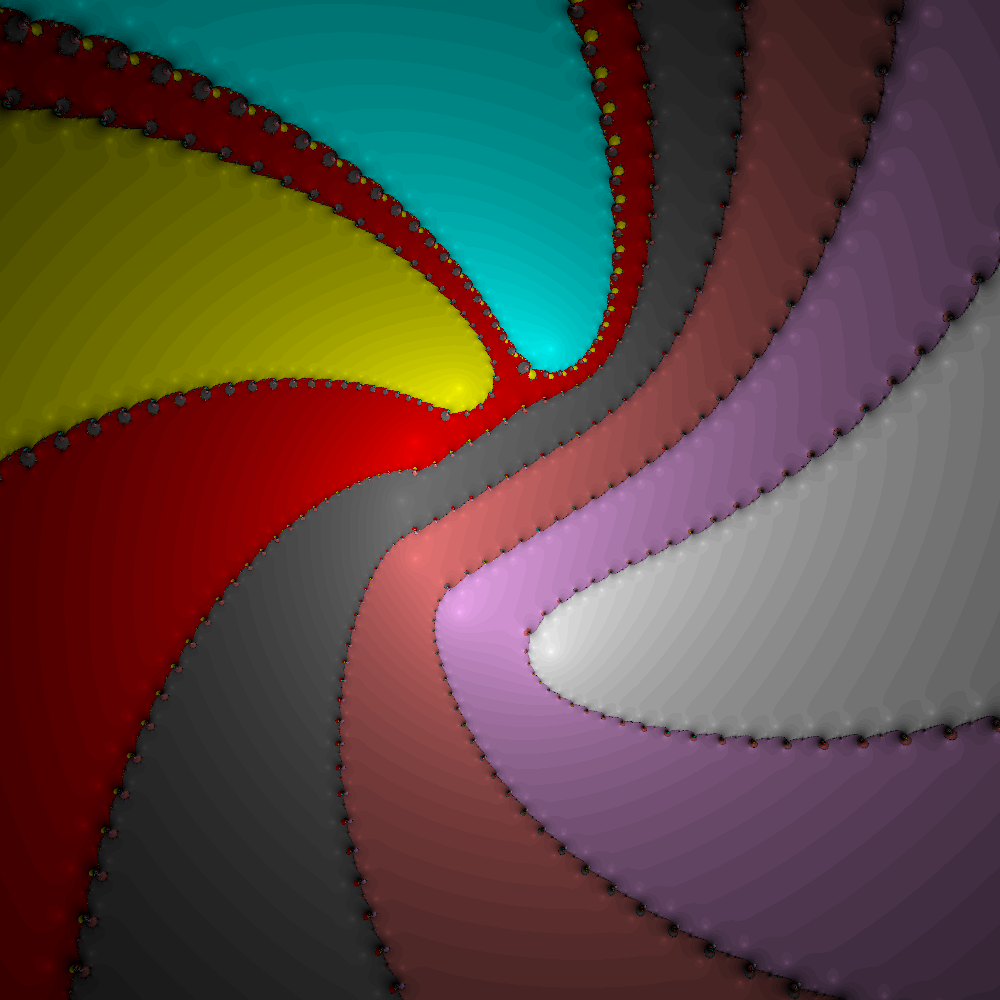}\

\vspace{0.1 cm}
\caption{Polynomiographs of $P_n(z)$, $n=2, \dots, 7$ under parameterized Newton.} \label{Fig3}
\end{center}
\end{figure}

Here we describe several polynomiographies for the first few partial sums based on the basic family. It is easy to show that the roots of $P_n(z)$ are simple, i.e. $P_n(z)$ and its derivative $P'_n(z)$ have no common zeros.  Also, it can be shown that the modulus of any root $\theta$ a root of $P_n(z)$ satisfies $0 <  |\theta | < n$.  The polynomiography of $P_1(z)=1+z$ is quite simple.  Any member of the basic family will converge to the root in one iteration. Polynomiographs of $P_n(z)$ under Newton's method for $n=2, \dots, 10$ are depicted in Figure \ref{Fig1},  showing them in increasing order from left to right and top to bottom.  The shape of the zeros form a convex shape,  reminiscent of  a parabola of the form $x=y^2$. This can be seen in the polynomiographs.  The norm of the roots goes to infinity as $n$ does. Figure \ref{Fig2} shows the polynomiography of $P_n(z)$, $n=2, \dots, 7$ under the point-wise convergence.  These are not fractal images.  Figure \ref{Fig3} shows the polynomiography of $P_n(z)$, $n=2, \dots, 7$ under parametrized Newton's method all for a particular value of $\alpha$.

\section*{Polynomiographs of Szeg\"o Partial Sums}
\begin{figure}[h!tbp]
\begin{center}
\includegraphics[width=1.5 in]{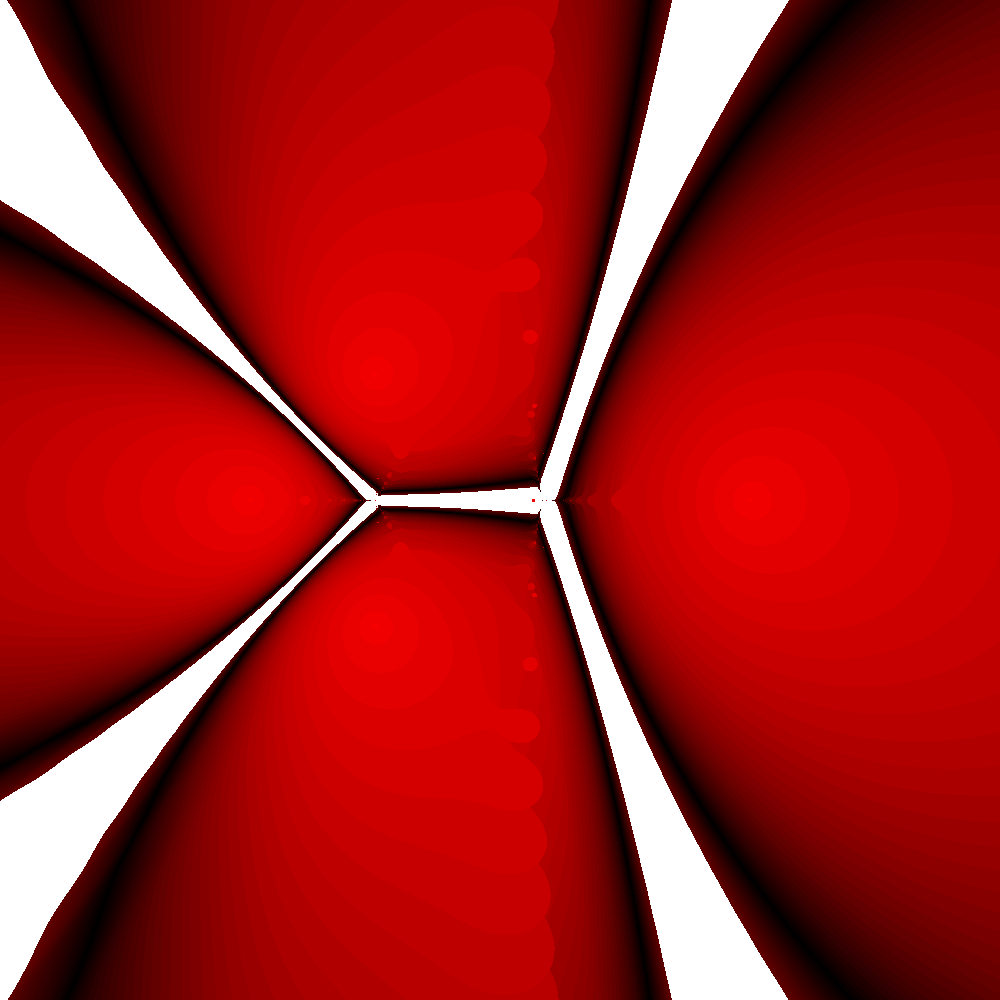} \includegraphics[width=1.5 in] {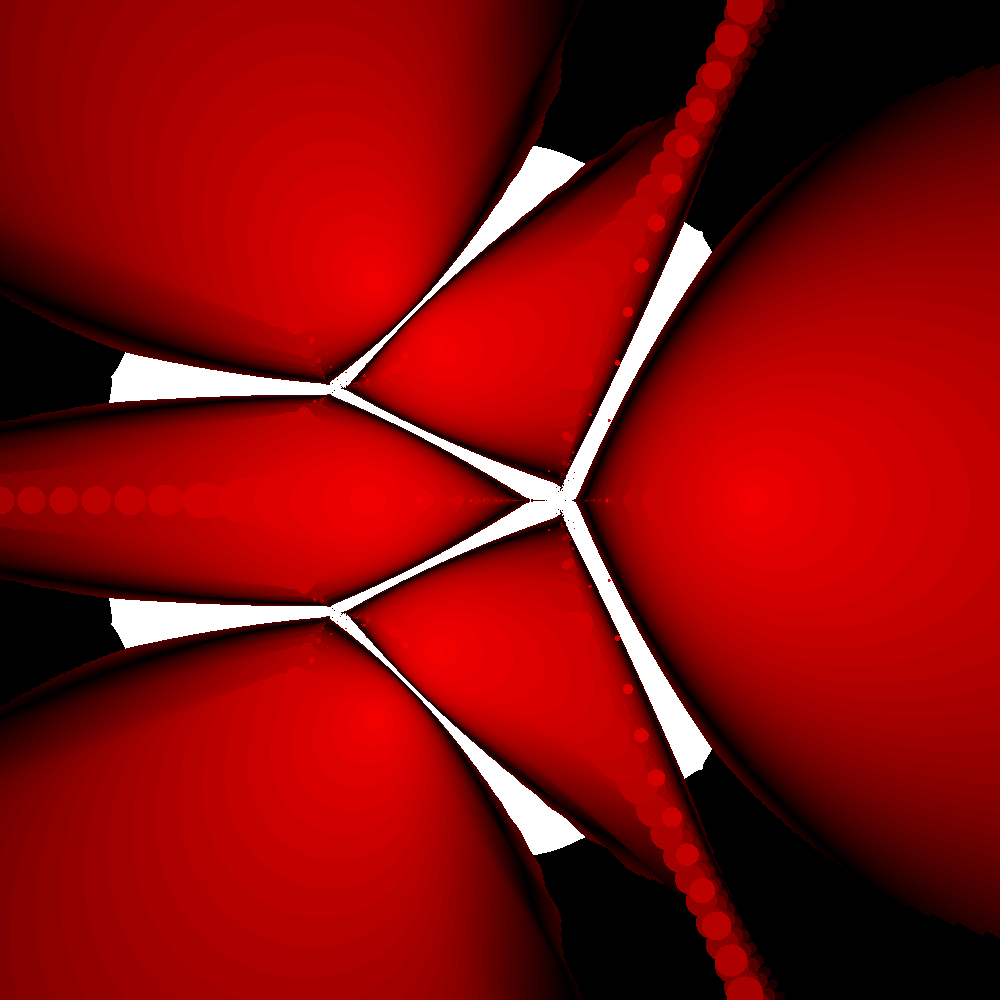} \includegraphics [width=1.5 in]{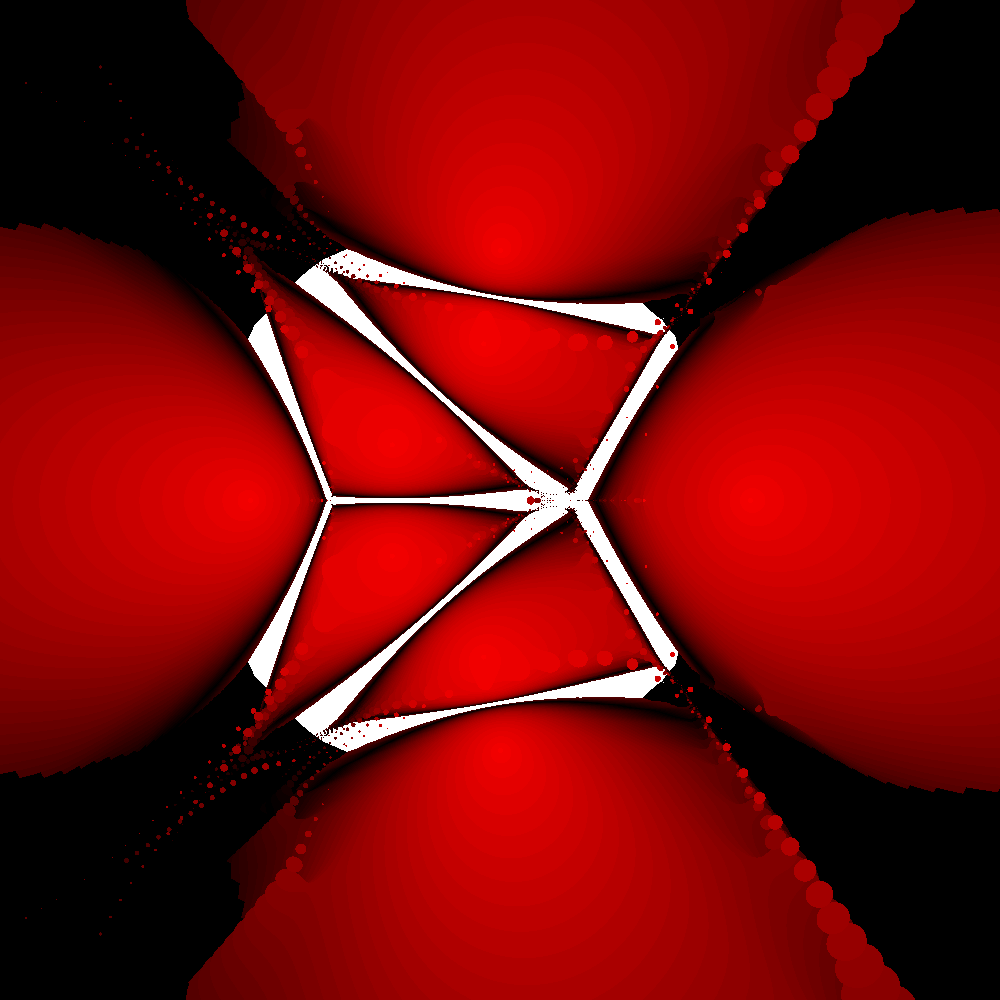}\

\vspace{0.1 cm}
\includegraphics[width=1.5 in]{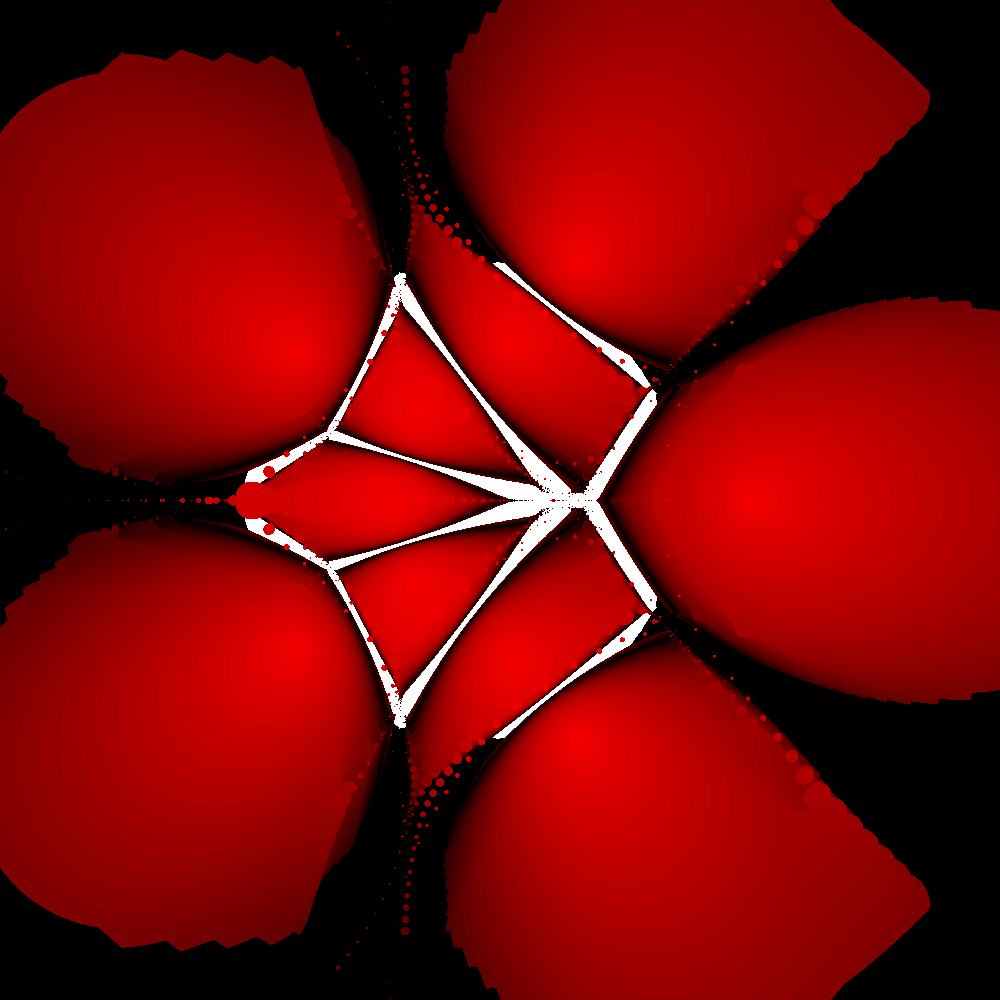} \includegraphics[width=1.5 in]{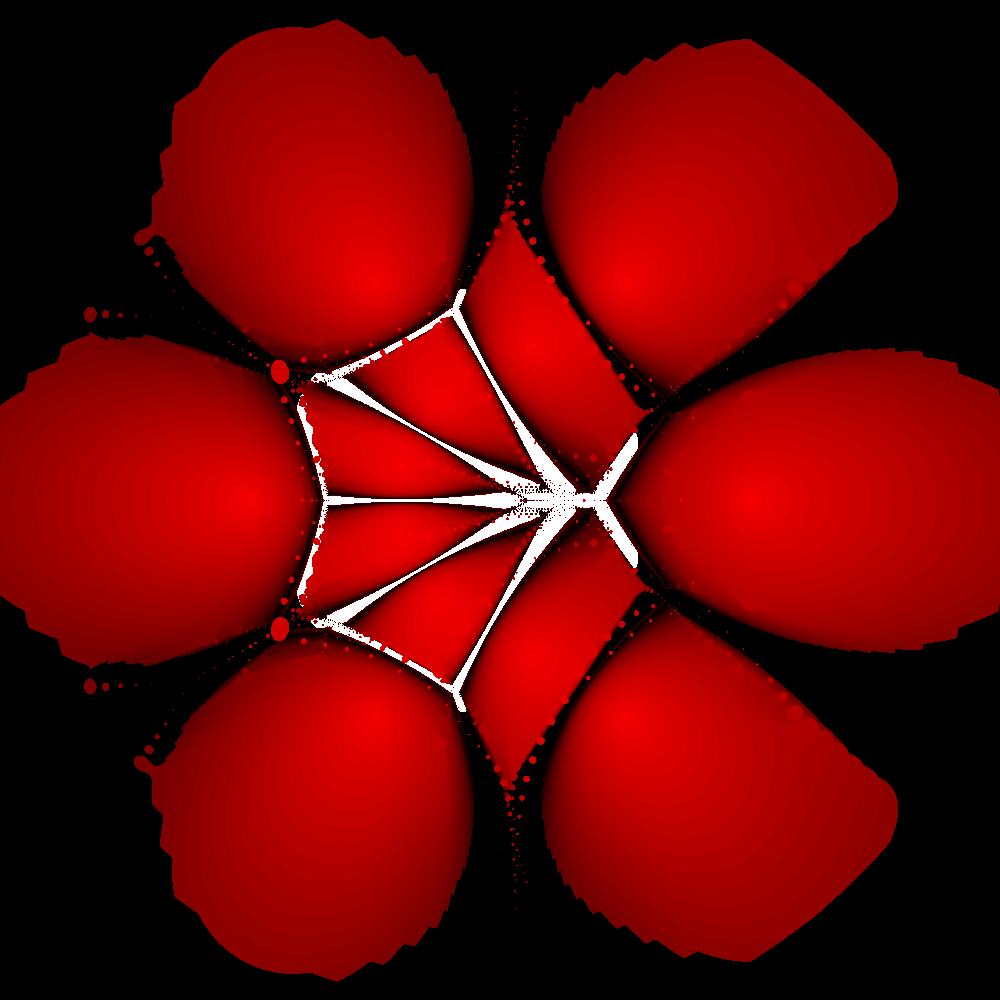} \includegraphics[width=1.5 in]{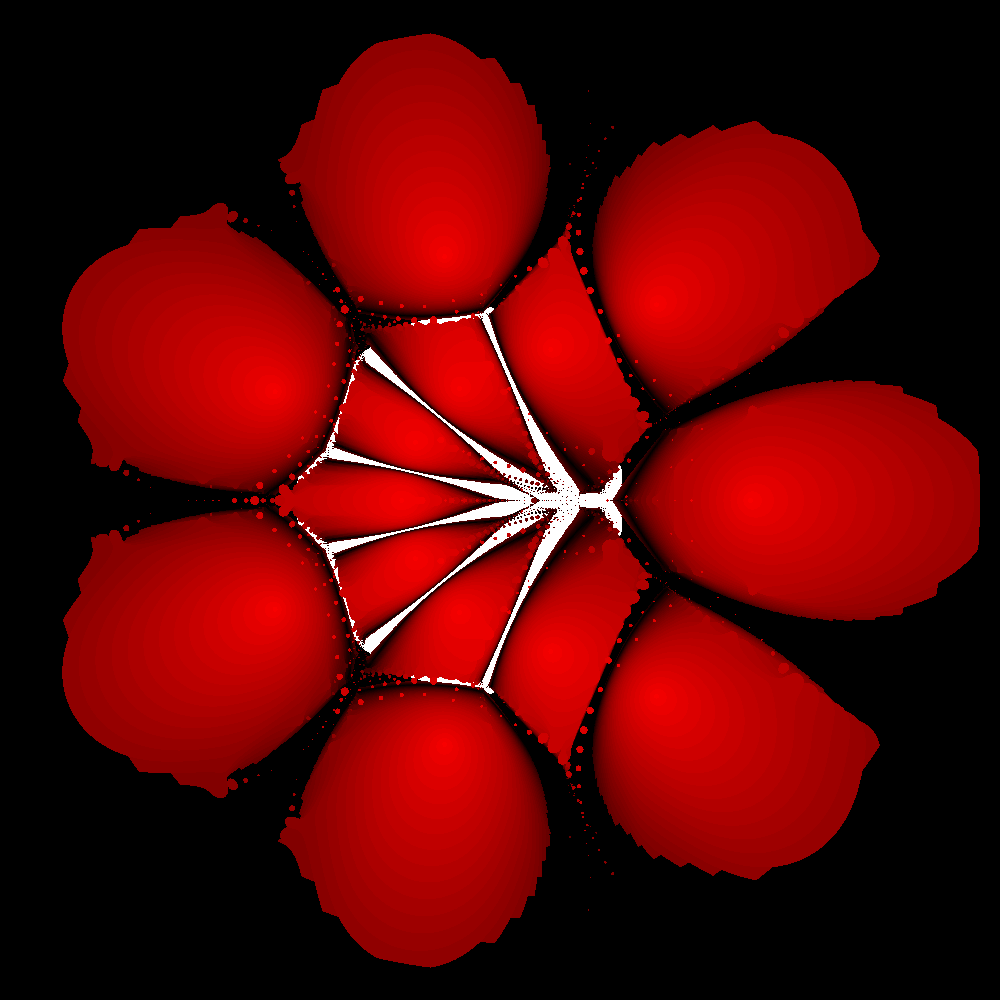}\

\vspace{0.1 cm}
\caption{Polynomiographs of $S_n(z) \times (z^n-1)$, $n=2, \dots, 7$ under point-wise convergence.} \label{Fig4}
\end{center}
\end{figure}

The norm of the roots of $P_n(z)$ get large as $n$ does. Szeg\"o partial sums are
\begin{equation}
S_n(z)=P_n(nz)= \sum_{k=0}^n \frac{(nz)^k}{k!}=1+nz+\frac{(nz)^2}{2!}+\cdots + \frac{(nz)^n}{n!}.
\end{equation}
Many interesting properties of this polynomial are known, see   P\'olya and G.~ Szeg\"o \cite{Pol} and Zemyan \cite{Zemyan}.
If $\theta$ is a root of $P_n(z)$ then $\theta/n$ is a root of $S_n(z)$. Thus all the roots of $S_n(z)$ are inside the disc of radius one, centered at the origin.  Polynomiographs of $S_n(z)$ for small $n$ look like scaled version of those of $P_n(z)$. However, as $n$ goes to infinity the zeros of $S_n(z)$ bend, forming an almond-shape inside the unit disc, see Zemyan \cite{Zemyan}.
Figure \ref{Fig4} shows the polynomiography of $S_n(z) \times (z^n-1)$ under point-wise convergence of the basic family.  The reason for multiplying by the roots of unity $z^n-1$ is two-fold: To show roots  lie inside the unit disc and  that under multiplication of polynomials we can generate interesting polynomiographs as science and as art.

\section*{Concluding Remarks}

In this article I have demonstrated polynomiography for the partial sums of the exponential series.  One can appreciate the images as art, but also as a way to get interested in learning or teaching root-finding algorithms.  What is intriguing about polynomiography software is that in the course of generating images we can learn about the shape of the zeros and get introduced to many other concepts in math and related areas.  Polynomiography is a medium for STEM, a bridge to learning or teaching about different subject areas, and making artistic images by considering variations of polynomials,  root-finding algorithms, coloring techniques, and operations such as multiplication of polynomials, scaling their zeros, compositions and more. Indeed we could make many artistic images based on the partial sums alone.   Interested educators can get a student module, see \cite{Ander}, as well as a link to a free demo polynomiography software upon registration at  \url{http://www.comap.com/Free/VCTAL/}.  See also, Choate \cite{jon} for lesson plans and a short manual for the demo software.



\end{document}